\documentclass[10pt]{amsart}         
\newcommand{\C}{{\mathbb C}}
\renewcommand{\P}{{\mathbb P}}
\def\fbar{\overline f} 

\newtheorem{theorem}{Theorem}

\newtheorem{proposition}[theorem]{Proposition}

\begin{document}
\title[Rational Polynomials in two variables]{Nontrivial Rational Polynomials in
two variables have reducible fibres}
\author{Walter D. Neumann}
\address{Department of Mathematics\\The University of
Melbourne\\Parkville, Vic 3052\\Australia}
\email{neumann@maths.mu.oz.au}
\author{Paul Norbury}
\address{Department of Mathematics\\The University of
Melbourne\\Parkville, Vic 3052\\Australia}
\email{norbs@maths.mu.oz.au}
\thanks{This research is supported by the Australian Research Council.}

\maketitle

We shall call a polynomial map $f\colon\C^2\to\C$ a ``coordinate'' if
there is a $g$ such that $(f,g)\colon\C^2\to\C^2$ is a polynomial
automorphism.  Equivalently, by Abhyankar-Moh and Suzuki, $f$ has one
and therefore all fibres isomorphic to $\C$.  Following
\cite{miyanishi-sugie} we call a polynomial $f\colon\C^2\to\C$
``rational'' if the general fibres of $f$ (and hence all fibres of
$f$) are rational curves.  The following theorem, which says that a
rational polynomial map with irreducible fibres cannot be part of a
counterexample to the $2$-dimensional Jacobian Conjecture, has
appeared in the literature several times. It appears with an algebraic
proof in Razar \cite{razar}. It is Theorem 2.5 of Heitmann
\cite{heitmann} (as corrected in the Corrigendum),
and L\^e and Weber, who give a geometric proof in \cite{le-weber},
also cite the reference Friedland \cite{friedland}, which we have not
seen.
\begin{theorem}
If $f\colon\C^2\to \C$ is a  rational polynomial map with
irreducible fibres and is not a coordinate then $f$ has no
jacobian partner (i.e., no polynomial $g$ such that the jacobian of
$(f,g)$ is a non-zero constant).
\end{theorem}

In this note we prove the above theorem is empty:
\begin{theorem}
There is no $f$ satisfying the assumptions of the above theorem. That
is, a  rational $f$ with irreducible fibres is a coordinate.
\end{theorem}
\begin{proof}
This theorem is implicit in \cite{miyanishi-sugie}.
Suppose $f$ is  rational.  As in \cite{miyanishi-sugie},
\cite{le-weber}, \cite{neumann2}, 
etc., we consider a nonsingular compactification
$Y=\C^2\cup E$ of $\C^2$ such that $f$ extends to a holomorphic map
$\fbar\colon Y\to\P^1$. Then $E$ is a union of smooth rational curves
$E_1,\dots,E_n$ with normal crossings.  An $E_i$ is called
\emph{horizontal} if $\fbar|E_i$ is non-constant.  Let $\delta$ be the
number of horizontal curves.  Then we have 
$$\delta-1=\sum_{a\in\C}(r_a-1),$$
where $r_a$ is the number of irreducible components of $f^{-1}(a)$.
This is Lemma 1.6 of Miyanishi and Sugie \cite{miyanishi-sugie} who
attribute it to Saito \cite{saito} and
Lemma 4 of L\^e-Weber \cite{le-weber} who attribute it to
Kaliman \cite{kaliman}, corollary 2. 
The proof is simple arithmetic from the topological
observation that on the one hand the euler characteristic of $Y$ is
$n+2$ and on the other hand it is $4+\sum_{a\in\P^1}(\overline r_a-1)
$, where $\overline r_a$ is the number of components of
$\fbar^{-1}(a), a\in\P^1$.

By this formula,
if $f$ has irreducible fibres then there is just one
horizontal curve. Lemma 1.7 of
\cite{miyanishi-sugie} now says that $f$ is a 
coordinate. This also follows from the following proposition, which
implies that the generic fibres of $f$ have just one point at infinity
and are thus isomorphic to $\C$.
\end{proof}

\begin{proposition}
Let $f\colon\C^2\to\C$ be any polynomial map and $\fbar\colon Y\to\P^1$
an extension as above. Denote by $d$ the greatest common divisor of
the degrees of $\fbar$ on the horizontal curves of $Y$ and $D$ the sum
of these degrees. Then the general fibre of $f$ has $d$ components (so
$f=h\circ f_1$ for some polynomials $f_1\colon\C^2\to \C$ and
$h\colon\C\to\C$ with $\operatorname{degree}(h)=d$), each of which is
a compact curve with $D/d$ punctures.
\end{proposition}

\begin{proof}
Let $E_1,\dots,E_\delta$ be the horizontal curves and
$d_1,\dots,d_\delta$ be the degrees of $\fbar$ on these.  Note that
the points at infinity of a general fibre $f^{-1}(a)$ are the points
where $\fbar^{-1}(a)$ meet the horizontal curves $E_i$, so there are
$d_i$ such points on $E_i$ for $i=1,\dots,\delta$. The relationship
between plumbing diagram and splice diagram (cf.\ \cite{neumann2,
eisenbud-neumann} says that the splice diagram $\Gamma$ for a regular
link at infinity for $f$ (cf.\ \cite{neumann1}) has $\delta$ nodes
with arrows at them, and the number of arrows at these nodes are
$d_1,\dots,d_\delta$ respectively.  Let $\Gamma_0$ be the same splice
diagram but with $d_1/d,\dots,d_\delta/d$ arrows at these nodes. Then
a minimal Seifert surface $S$ for the link represented by $\Gamma$
will consist of $d$ parallel copies of a minimal Seifert surface for
the link represented by $\Gamma_0$, so this $S$ has $d$
components. But the general fibre of $f$ is such a minimal Seifert
surface (\cite{neumann1}, Theorem 1), completing the proof. (It also
follows that $\Gamma_0$ is the regular splice diagram for the
polynomial $f_1$ of the proposition.)
\end{proof}

\end{document}